\documentclass[12pt,]{amsart}
\usepackage{amsmath}
\usepackage{amscd}
\usepackage{amssymb}
\usepackage{latexsym}
\input xy
\xyoption{all} \CompileMatrices
\newtheorem{theorem}{Theorem}[section]
\newtheorem{lemma}[theorem]{Lemma}
\newtheorem{corollary}[theorem]{Corollary}
\newtheorem{proposition}[theorem]{Proposition}

\newtheorem{remark}[theorem]{Remark}
\newtheorem{definition}[theorem]{Definition}

\newcommand{\bgl}{\begin{equation}}         %eine Gleichung mit Ziffer
\newcommand{\egl}{\end{equation}}
\newcommand{\bgln}{\begin{eqnarray}}        %mehrere Gleichungen mit Ziffer
\newcommand{\egln}{\end{eqnarray}}
\newcommand{\bglnoz}{\begin{eqnarray*}}     %mehrere Gleichungen ohne Ziffer
\newcommand{\eglnoz}{\end{eqnarray*}}
\newcommand{\btheo}{\begin{theorem}}
\newcommand{\etheo}{\end{theorem}}
\newcommand{\blemma}{\begin{lemma}}
\newcommand{\elemma}{\end{lemma}}
\newcommand{\bproof}{\begin{proof}}
\newcommand{\eproof}{\end{proof}}
\newcommand{\bbew}{\begin{beweis}}
\newcommand{\ebew}{\end{beweis}}
\newcommand{\bremark}{\begin{remark}\em}
\newcommand{\eremark}{\end{remark}}
\newcommand{\bdefin}{\begin{definition}}
\newcommand{\edefin}{\end{definition}}
\newcommand{\bprop}{\begin{proposition}}
\newcommand{\eprop}{\end{proposition}}
\newcommand{\bcor}{\begin{corollary}}
\newcommand{\ecor}{\end{corollary}}
\newcommand{\mn}{\par\medskip\noindent}

\newcommand{\cC}{\mathcal C}
\newcommand{\cD}{\mathcal D}

\newcommand{\cF}{\mathcal F}

\newcommand{\cK}{\mathcal K}

\newcommand{\cM}{\mathcal M}

\newcommand{\cO}{\mathcal O}
\newcommand{\cP}{\mathcal P}
\newcommand{\cQ}{\mathcal Q}

\newcommand{\lori}{\longrightarrow}

\def\SEMI{\mbox{$\times\kern-2pt\vrule height5pt width.6pt \kern3pt $}}

\newcommand{\Img}{{\rm Im\,}}

\newcommand{\Ad}{{\rm Ad\,}}

\def\Az{\mathbb{A}}
\def\Cz{\mathbb{C}}

\def\Nz{\mathbb{N}}

\def\Qz{\mathbb{Q}}
\def\Rz{\mathbb{R}}
\def\Tz{\mathbb{T}}
\def\Zz{\mathbb{Z}}

\title{$C^*$-algebras associated with the $ax+b$-semigroup over $\Nz$}
\author{Joachim Cuntz}
\address{Joachim Cuntz, Mathematisches Institut, Einsteinstr.62, 48149 M\"unster, Germany}
\email{cuntz@math.uni-muenster.de}
\urladdr{http://www.math.uni-muenster.de/u/cuntz/cuntz}
\thanks{Research supported by the Deutsche Forschungsgemeinschaft}
\subjclass[2000]{Primary: 58B34, 46L05} \keywords{purely
infinite, Bost-Connes algebra}

\begin{document}

\maketitle \begin{abstract}\noindent We present a
$C^*$-algebra which is naturally associated to the
$ax+b$-semigroup over $\mathbb N$. It is simple and purely
infinite and can be obtained from the algebra considered by
Bost and Connes by adding one unitary generator which
corresponds to addition. Its stabilization can be described
as a crossed product of the algebra of continuous functions,
vanishing at infinity, on the space of finite adeles for
$\mathbb Q$ by the natural action of the $ax+b$-group over
$\mathbb Q$.
\end{abstract}\tableofcontents
\newpage
\section{Introduction}
In this note we present a $C^*$-algebra (denoted by
$\cQ_\Nz$) which is associated to the $ax+b$-semigroup over
$\Nz$. It is in fact a natural quotient of the $C^*$-algebra
associated with this semigroup by some additional relations
(which make it simple and purely infinite). These relations
are satisfied in representations related to number theory.
The study of this algebra is motivated by the construction of
Bost-Connes in \cite{BoCo}. Our $C^*$-algebra contains the
algebra considered by Bost-Connes, but in addition a
generator corresponding to translation by the additive group
$\Zz$.

As a $C^*$-algebra, $\cQ_\Nz$ has an interesting structure.
It is a crossed product of the Bunce-Deddens algebra
associated to $\Qz$ by the action of the multiplicative
semigroup $\Nz^\times$. It has a unique canonical KMS-state.
We also determine its $K$-theory, whose generators turn out to
be determined by prime numbers.

On the other hand, $\cQ_\Nz$ can also be obtained as a crossed
product of the commutative algebra of continuous functions on
the completion $\hat{\Zz}$ by the natural action of the
$ax+b$-semigroup over $\Nz$. More interestingly, its
stabilization is isomorphic to the crossed product of the
algebra $\cC_0(\Az_f)$ of continuous functions on the space
of finite adeles by the natural action of the $ax+b$-group
$P^+_\Qz$ over $\Qz$. Somewhat surprisingly one obtains
exactly the same $C^*$-algebra $\cQ_\Nz$ (up to
stabilization) working with the completion $\Rz$ of $\Qz$ at
the infinite place and taking the crossed product by the
natural action of the $ax+b$-group $P^+_\Qz$ on $\cC_0(\Rz)$.

In the last section we consider the analogous construction of
a $C^*$-algebra replacing the multiplicative semigroup
$\Nz^\times$ by $\Zz^\times$, i.e. omitting the condition of
positivity on the multiplicative part of the
$ax+b$-(semi)group. We obtain a purely infinite $C^*$-algebra
$\cQ_\Zz$ which can be written as a crossed product of
$\cQ_\Nz$ by $\Zz/2$. The fixed point algebra for its
canonical one-parameter group $(\lambda_t)$ is a dihedral
group analogue of the Bunce-Deddens algebra. Its $K$-theory
involves a shift of parity from $K_0$ to $K_1$ and vice
versa. Its stabilization is isomorphic to the natural crossed
product $\cC_0(\Az_f)\rtimes P_\Qz$.

\section{A canonical representation of the
$ax+b$-semigroup over $\Nz$}\label{canrep} We denote by $\Nz$ the
set of natural numbers including 0. $\Nz$ will normally be regarded
as a semigroup with addition.

We denote by $\Nz^\times$ the set of natural numbers
excluding 0. $\Nz^\times$ will normally be regarded as a
semigroup with multiplication.

The natural analogue, for a semigroup $S$, of an unitary
representation of a group is a representation of $S$ by
isometries, i.e. by operators $s_g, g\in S$ on a Hilbert
space that satisfy $s_g^*s_g=1$.

On the Hilbert space $\ell^2(\Nz)$ consider the isometries
$s_n, n\in \Nz^\times$ and $v^k, k \in \Nz$ defined by

$$ v^k(\xi_m)=\xi_{m+k}\qquad s_n (\xi_m) =\xi_{mn}$$

where $\xi_m, m\in \Nz$ denotes the standard orthonormal
basis.

We have $s_ns_m =s_{nm}$ and $v^nv^m=v^{n+m}$, i.e. $s$ and
$v$ define representations of $\Nz^\times$ and $\Nz$
respectively, by isometries. Moreover we have the following
relation
$$s_nv^k = v^{nk}s_n $$
which expresses the compatibility between multiplication and
addition.\mn In other words the $s_n$ and $v^k$ define a
representation of the $ax+b$-semigroup
$$P_{\Nz} = \{ \left(\begin{array}{cc}1&k\\0&n\end{array}\right)
\mathop{|} n\in\Nz^\times,k\in\Nz\}$$ over $\Nz$, where the
matrix $\left(\begin{array}{cc}1&k\\0&n\end{array}\right)$ is
represented by $v^k s_n$.

Note that, for each $n$, the operators $s_n, vs_n, \ldots,
v^{n-1}s_n$ generate a $C^*$-algebra isomorphic to $\cO_n$,
\cite{CuCMP}.

\section{A purely infinite simple $C^*$-algebra associated with the
$ax+b$-semigroup} The $C^*$-algebra $A$ generated by the
elements $s_n$ and $v$ considered in section \ref{canrep}
contains the algebra $\cK$ of compact operators on
$\ell^2\Nz$. Denote by $u$, resp. $u^k$ the image of $v$,
resp. $v^k$ in the quotient $A/\cK$. We also still denote by
$s_n$ the image of $s_n$ in the quotient. Then the $u^k$ are
unitary and are defined also for $k\in \Zz$ and they
furthermore satisfy the characteristic relation
$\sum_{k=0}^{n-1}u^ke_nu^{-k}=1$ where $e_n=s_ns_n^*$ denotes
the range projection of $s_n$. This relation expresses the
fact that $\Nz$ is the union of the sets of numbers which are
congruent to $k$ mod $n$ for $k=0,\ldots,n-1$.

We now consider the universal $C^*$-algebra generated by
elements satisfying these relations.

\bdefin We define the $C^*$-algebra $\cQ_\Nz$ as the
universal $C^*$-algebra generated by isometries $s_n, n\in
\Nz^\times$ with range projections $e_n=s_ns_n^*$, and by a
unitary $u$ satisfying the relations
$$s_ns_m=s_{nm}, \qquad s_nu=u^ns_n,\qquad \sum_{k=0}^{n-1}u^ke_nu^{-k}=1$$ for $n,m \in
\Nz^\times$.\edefin

\blemma\label{comm} In $\cQ_\Nz$ we have
\begin{itemize}
\item[(a)] $e_n =\sum_{i=0}^{m-1} u^{in}e_{nm}u^{-in}$ for
all $n,m\in \Nz^\times$.
\item[(b)] $e_ps_q =e_{pq}s_p=s_qe_p$ and
$e_pe_q=e_{pq}=e_qe_p$ when $p$ and $q$ are relatively prime.
\item[(c)] $s_n^*s_m =s_ms_n^*$ for all $n,m$.
\end{itemize} \elemma
\bproof (a) This follows by conjugating the identity $1 =
\sum u^ie_mu^{-i}$ by $s_n\sqcup s_n^*$ and using the fact
that $s_ne_ms_n^*=e_{nm}$.\\
(b) Since the $u^ie_{pq}u^{-i}$ are pairwise orthogonal for
$0\leq i<pq$, we see that $u^{lp}e_{pq}u^{-lp}\perp
u^{kq}e_{pq}u^{-kq}$ if $0< l<q, 0< k<p$. Thus, using (a),
$$e_ps_q= \sum_{l=0}^{q-1}u^{lp}e_{pq}u^{-lp}\sum_{k=0}^{p-1}u^{kq}e_{pq}u^{-kq}
s_q = e_{pq}s_q$$
This obviously implies $e_pe_q=e_{pq}$ and, by symmetry
$e_qe_p=e_{pq}$. In particular, $e_p$ and $e_q$ commute.\\
(c) Using (b) we get
$$s_p(s_p^*s_q)s_q^*=e_pe_q=e_{pq}=s_{pq}s_{pq}^*=s_p(s_qs_p^*)s_q^*$$
Since $s_p, s_q$ are isometries this implies that the
expressions in parentheses on the left and right hand side
are equal. Here we have, in a first step, assumed that $p,q$
are prime. However any $s_n$ is a product of $s_p$'s with $p$
prime. \eproof

From \ref{comm} (a) it follows that any two of the projections
$u^ie_nu^{-i}$ commute. We denote by $D$ the commutative
subalgebra of $\cQ_\Nz$ generated by all these projections.

We also denote by $\cF$ the subalgebra of $\cQ_\Nz$ generated
by $u$ and the projections $e_n$, $n\in \Nz^\times$.

To analyze the structure of $\cQ_\Nz$ further we write it as
an inductive limit of the subalgebras $B_n$ generated by
$s_{p_1}, s_{p_2}, \ldots , s_{p_n}$ and $u$, where $p_1,p_2,
\ldots , p_n$ denote the $n$ first prime numbers. Each $B_n$
contains a natural (maximal) commutative subalgebra $D_n$
generated by all projections of the form $u^ke_mu^{*k}$ where
$m$ is a product of powers of the $p_1,\ldots, p_n$ (i.e. a
natural number that contains only the $p_i$ as prime
factors). \blemma\label{spec} The spectrum {\rm Spec}$D_n$ of
$D_n$ can be identified canonically with the compact space
$$ \{0,\ldots , p_1-1\}^\Nz\times\ldots \times \{0,\ldots ,
p_n-1\}^\Nz\cong \hat{\Zz}_{p_1}\times\ldots\times
\hat{\Zz}_{p_n}$$ \elemma \bproof $D_n$ is the inductive
limit of the subalgebras $D_n^{(k)}\cong \Cz^{l_k}$ with
$l_k=p_1^kp_2^k\ldots p_n^k$. The algebra $D_n^{(k)}$ is
generated by the pairwise orthogonal projections
$u^ie_{l_k}u^{-i}, 0\leq i< l_k$ and in fact, by Lemma
\ref{comm} (a), $D_n^{(k)}$ is the $k$-fold tensor product of
$D_n^{(1)}$ by itself. \eproof

Consider the action of $\Tz^n$ on $B_n$ given by
$$\alpha_{(t_1,\ldots,t_n)}(s_{p_i})=t_is_{p_i}$$ and denote
by $\cF_n$ the fixed-point algebra for $\alpha$ (i.e.
$\cF_n=\cF\cap B_n$). There is a natural faithful conditional
expectation $E: B_n\to \cF_n$ defined by $E(x)=
\int_{\Tz^n}\alpha_t(x)dt$.

Now $D_n$ is the fixed point algebra, in $\cF_n$, for the
action $\beta$ of $\Tz$ on $\cF_n$ given by
$$\beta_t (e_k)=e_k\qquad \beta_t(u)=e^{it}u$$
and there is an associated expectation $F:\cF_n\to D_n$
defined by $F(x)= \int_{\Tz}\beta_t(x)dt$. The composition
$G=F\circ E$ gives a faithful conditional expectation $A_n\to
D_n$. These conditional expectations extend to the inductive
limit and thus give conditional expectations $E:\cQ_\Nz\to
\cF$, $F:\cF\to \cD$ and $G:\cQ_\Nz \to \cD$.

\btheo\label{pi} The $C^*$-algebra $\cQ_\Nz$ is simple and
purely infinite.\etheo

\bproof Since inductive limits of purely infinite simple
$C^*$-algebras are purely infinite simple \cite{Ror}, 4.1.8
(ii), it suffices to show that each $B_n$ is purely infinite
simple.

For each $N$, denote, as above, by $D_n^{(N)}$ the subalgebra
of $D_n$ generated by $\{u^ke_lu^{-k}, k\in \Zz\}$ where
$l=p_1^Np_2^N\ldots p_n^N$. The natural map Spec$ D_n\to$
Spec $D_n^{(N)}$ is surjective and, by the proof of Lemma
\ref{spec}, $D_n^{(N)}\cong \Cz^l$.

Choose $\xi_1,\xi_2,\ldots,\xi_l$ in Spec $D_n$ such that
$\{\xi_1,\xi_2,\ldots,\xi_l\}\to$ Spec $D_n^{(N)}$ is
bijective and such that $$\xi_i (s_k \sqcup s_k^*) \neq
\xi_i(\sqcup)$$ for all $k$ of the form $k=p_1^{m_1}\ldots
p_n^{m_n}$ with $m_i\leq N$. We can choose pairwise
orthogonal projections $f_1,f_2,\ldots,f_l$ in $D_n$ with
sufficiently small support around the $\xi_i$ such that
$f_is_kf_i=0$ for all $1\leq i\leq l$ and $k$ as above and
such that $f_ixf_i=\xi_i(x)f_i$ for all $x\in D_n^N$. Then
the map $\varphi :D_n^{(N)}\to C^*(f_1,\ldots,f_l)\cong
\Cz^l$ defined by $x\mapsto \sum f_i xf_i$ is an
isomorphism.\mn

Denote by $\cP^{(N)}$ the set of linear combinations of all
products of the form $u^ks_{p_1}^{m_1}\ldots s_{p_n}^{m_n}$
with $m_i\leq N$. For $x\in \cP^{(N)}$ we have that
$$\varphi(G(x))=\sum f_i xf_i$$
Let now $0\leq x\in B_n$ be different from 0. Since $G$ is
faithful, $G(x)\neq 0$ and we normalize $x$ such that $\|
G(x)\|=1$. Let $y\in \cP^{(N)}$, for sufficiently large $N$,
be such that $\|x-y\|<\varepsilon\leq 1$. We may also assume
that $\| G(y)\|=1$. Then there exists $i_0, 1\leq i_0\leq l$
such that $f_{i_0}yf_{i_0}=f_{i_0}$. Moreover, there exists
an isometry $s\in B_n$ (of the form $s=u^ks_m$) such that
$s^*f_{i_0}s=1$. \mn In conclusion, we have $s^*ys=1$ and
$$\|s^*xs-1\|= \|s^*xs-s^*ys\| = \|s^*(x-y)s\| < \varepsilon$$
This shows that $s^*xs$ is invertible. \eproof

As a consequence of the simplicity of $\cQ_\Nz$ we see that
the canonical representation on $\ell^2(\Zz)$ by
$$s_n(\xi_k)=\xi_{nk}\qquad u^n(\xi_k)=\xi_{k+n}$$ is faithful.
Similarly, if we divide the $C^*$-algebra generated by the
analogous isometries on $\ell^2(\Nz)$, discussed in section
\ref{canrep}, by the canonical ideal $\cK$, we get an algebra
isomorphic to $\cQ_\Nz$.

The subalgebra $\cF$ is generated by $u$ and the projections
$e_n$, thus by a weighted shift. Thus we recognize $\cF$ as
the Bunce-Deddens algebra of the type where every prime
appears with infinite multiplicity. This well known algebra
has been introduced in \cite{BuDe} in exactly that form. It is
simple and has a unique tracial state. It can also be
represented as inductive limit of the inductive system $(M_n
\cC (S^1))$ with maps
$$M_n(\cC (S^1)) \lori M_{nk}(\cC (S^1))$$
mapping the unitary $u$ generating $\cC (S^1)$ to the $k\times
k$-matrix
$$\left(\begin{array}{ccccc}
0&0&\ldots&&u\\
1&0&\ldots&&0\\
0&1&\ldots&&0\\
&&\ldots&&\\
0&\ldots&&1&0 \end{array}\right)$$ From this description it
immediately follows that $K_0(\cF)=\Qz$ and $K_1(\cF)=\Zz$.

\bremark Just as $\cO_n$ is a crossed product of a
$UHF$-algebra by $\Nz$, we see that $\cQ_\Nz$ is a crossed
product $\cF\rtimes \Nz^\times$ by the multiplicative
semigroup $\Nz^\times$. The algebra $\cQ_\Nz$ also contains
the commutative subalgebra $D$. Since $D$ is the inductive
limit of the $D_n$, we see from \ref{spec} that Spec
$D=\hat{\Zz}:=\prod_p\hat{\Zz}_p$. The Bost-Connes algebra
$C_\Qz$ \cite{BoCo} can be described as a crossed product
$\cC(\hat{\Zz})\rtimes \Nz^\times$ and is a natural
subalgebra of $\cQ_\Nz$ (using the natural inclusion
$\cC(\hat{\Zz})\to \cF$).

We can also obtain $\cQ_\Nz$ by adding to the generators
$\mu_n, n\in \Nz$ (our $s_n$) and $e_\gamma , \gamma \in
\Qz/\Zz$ for $C_\Qz$, described in \cite{BoCo}, Proposition
18, one additional unitary generator $u$ satisfying
$$ue_\gamma = \gamma e_\gamma u \qquad u^n\mu_n = \mu_n u$$
(here we identify an element $\gamma$ of $\Qz/\Zz$ with the
corresponding complex number of modulus 1). \eremark

\section{The canonical action of
$\Rz$ on $\cQ_\Nz$} So far, in our discussion, the isometries
associated with each prime number appear as generators on the
same footing and it is a priori not clear how to determine,
for two prime numbers $p$ and $q$, the size of $p$ and $q$
(or even the bigger number among $p$ and $q$) from the
corresponding generators $s_p$ and $s_q$. In fact, the
$C^*$-algebra generated by the $s_n$, $n\in \Nz$ is the
infinite tensor product of one Toeplitz algebra for each
prime number and in this $C^*$-algebra there is no way to
distinguish the $s_p$ for different $p$.

However, the fact that we have added $u$ to the generators
allows to retrieve the $n$ from $s_n$ using the KMS-condition.

\bdefin Let $(A,\lambda_t)$ be a $C^*$-algebra equipped with
a one-parameter automorphism group $(\lambda_t)$, $\tau$ a
state on $A$ and $\beta\in (0,\infty]$. We say that $\tau$
satisfies the $\beta$-KMS-condition with respect to
$(\lambda_t)$, if for each pair $x, y$ of elements in $A$,
there is a holomorphic function $F_{x,y}$, continuous on the
boundary, on the strip $\{z\in \Cz \mathop{|} \Img z \in
[0,\beta]\}$ such that $$F_{x,y}(t) =\tau (x\lambda_t (y))
\qquad F_{x,y}(t+i\beta) =\tau (\lambda_t (y)x)$$ for $t\in
\Rz$.\edefin

\bprop Let $\tau_0$ be the unique tracial state on $\cF$ and
define $\tau$ on $\cQ_\Nz$ by $\tau = \tau_0\circ E$. Let
$(\lambda_t)$ denote the one-parameter automorphism group on
$\cQ_\Nz$ defined by $\lambda_t(s_n)= n^{it}s_n$ and
$\lambda_t (u)=u$. Then $\tau$ is a 1-KMS-state for
$(\lambda_t)$.\eprop \bproof According to \cite{BrRo} 5.3.1,
it suffices to check that $$\tau (x\lambda_{i}(y))=\tau (yx)$$
for a dense *-subalgebra of analytic vectors for
$(\lambda_t)$. Here $(\lambda_z)$ denotes the extension to
complex variables $z\in \Cz$ of $(\lambda_t)$ on the set of
analytic vectors.\mn However it is immediately clear that
this identity holds for $x, y$ linear combinations of
elements of the form $as_n$ or $s_m^*b$, $a,b\in \cF$. Such
linear combinations are analytic and dense. \eproof

\btheo There is a unique state $\tau$ on $\cQ_\Nz$ with the
following property:\\there exists a one-parameter
automorphism group $(\lambda_t)_{t\in \Rz}$ for which $\tau$
is a 1-KMS-state and such that $\lambda_t (u)=u$, $\lambda_t
(e_n)=e_n$ for all $n$ and $t$. Moreover we have
\begin{itemize}
\item $\tau$ is given by $\tau=\tau_0\circ E$ where $\tau_0$
is the canonical trace on $\cF$.
\item the one-parameter group for which $\tau$ is 1-KMS,
is unique and is the standard automorphism group considered
above, determined by
$$\lambda_t (s_n) = n^{it}s_n\qquad \lambda_t (u) = u$$
\end{itemize} \etheo
\bproof Since $\lambda_t$ acts as the identity on $e_n$ and
on $u$, it is the identity on $\cF =C^*(u, \{e_n\})$. If
$\tau$ is a KMS-state for $(\lambda_t)$, it therefore has to
be a trace on $\cF$. It is well-known (and clear) that there
is a unique trace $\tau_0$ on $\cF$. For instance, the
relation
$$\sum_{0\leq k<n}u^k e_n u^{-k}=1$$
shows that $\tau_0(e_n)=1/n$ and $\tau_0(u)=0$.

The relation
$$\tau (s_n^*s_n) = n \tau (s_ns_n^*)$$ and the 1-KMS-condition show
that $\lambda_t (s_n) = n^{it}s_n$. \eproof

\section{The $K$-groups of $\cQ _\Nz$}
The $K$-groups of $\cQ_\Nz$ can be computed using the fact
that $\cQ _\Nz$ is a crossed product of the Bunce-Deddens
algebra $\cF$ by the semigroup $\Nz^\times$. The $K$-groups
of a Bunce-Deddens algebra are well known and easy to
determine using its representation as an inductive limit of
algebras of type $M_n(\cC(S^1))$. Specifically, for $\cF$, we
have
$$ K_0(\cF)= \Qz\qquad K_1(\cF)= \Zz$$

We can consider $\cQ_\Nz$ as an inductive limit of the
subalgebras $B_n =C^*(\cF, s_{p_1},\ldots,s_{p_n})$ where
$2=p_1<p_2<\ldots$ is the sequence of prime numbers in
natural order. Now $B_{n+1}$ can be considered as a crossed
product of $B_n$ by the action of the semigroup $\Nz$ given by
conjugation by $s_{p_{n+1}}$ (in fact $B_{n+1}$ is Morita
equivalent to a crossed product by $\Zz$ of an algebra Morita
equivalent to $B_n$), just as in \cite{CuCMP}. Thus the
$K$-groups of $B_n$ can be determined inductively using the
Pimsner-Voiculescu sequence.

\btheo The $K$-groups of $B_n$ are given by
$$K_0(B_n)\cong \Zz^{2^{n -1}}\qquad K_1(B_n)\cong \Zz^{2^{n
-1}}$$ \etheo \bproof The first application of the
Pimsner-Voiculescu sequence gives an exact sequence
$$\lori K_1(B_1) \lori \Qz \mathop{\lori}\limits^{id - \alpha_*}  \Qz \lori K_0(B_1)
\lori \Zz \mathop{\lori}\limits^{id - \alpha_*}
 \Zz \lori K_1(B_1) \lori$$

where $\alpha_*$ is the map induced by Ad $s_1$ on $K_0, K_1$.

Since $\alpha_*=2$ on $K_0(\cF)$ and $\alpha_*= 1$ on
$K_1(\cF)$ this gives $K_0(B_1)=\Zz$, $K_1(B_1)=\Zz$. In the
following steps $\alpha_i=$ Ad $s_{p_i}$ induces 1 on $K_0$
and on $K_1$. Thus each consecutive prime $p_i$ doubles the
number of generators of $K_0$ and $K_1$. \eproof

\section{Representations as crossed products}
For each $n$, define the endomorphism $\varphi_n$ of $\cQ_\Nz$
by $\varphi_n(x)= s_n xs_n^*$. Since
$\varphi_n\varphi_m=\varphi_{nm}$ this defines an inductive
system.

\bdefin We define $\overline{\cQ}_\Nz$ as the inductive limit
of the inductive system $(\cQ_\Nz, \varphi_n)$.\edefin

By construction we have a family $\iota_n$ of natural
inclusions of $\cQ_\Nz$ into the inductive limit
$\overline{\cQ}_\Nz$ satisfying the relations
$\iota_{nm}\varphi_n=\iota_m$. We denote by $1_k$ the element
$\iota_k(1)$ of $\overline{\cQ}_\Nz$. We have that
$1_k=\iota_{kl}(e_l)\leq e_{kl}$. The union of the
subalgebras $1_k\overline{\cQ}_\Nz1_k$ is dense in
$\overline{\cQ}_\Nz$ and $1_k\leq 1_{kl}$ for all $k,l$. In
order to define a multiplier $a$ of $\overline{\cQ}_\Nz$ it
therefore suffices to define $a1_k$ and $1_ka$ for all $k$.

We can extend the isometries $s_n$ naturally to unitaries
$\bar{s}_n$ in the multiplier algebra
$\cM(\overline{\cQ}_\Nz)$ of $\overline{\cQ}_\Nz$ by requiring
$$\bar{s}_n1_k=\iota_k(s_n)\qquad 1_k \bar{s}_n = \iota_{kn}(e_ks_n)$$
Note that this is well defined because, using \ref{comm}, we
have
$$
(1_k\bar{s}_n)1_l=\iota_{kl}(e_ks_n)\iota_{kl}(e_l)=\iota_{kl}(s_ne_ke_l)=
\iota_{kl}(e_ls_ks_ns_k^*)=\iota_k(1)\iota_l(s_n)=1_k(\bar{s}_n1_l)$$

\bprop The elements $\bar{s}_n$, $n\in \Nz$, define unitaries
in $\cM (\overline{\cQ}_\Nz)$ such that
$\bar{s}_n\bar{s}_m=\bar{s}_{nm}$ and such that
$\bar{s}_n\bar{s}_m^*=\bar{s}_m^*\bar{s}_n$. \\
Defining $\bar{s}_a=\bar{s}_n\bar{s}_m^*$ for $a=n/m\in
\Qz^\times_+$ we define unitaries in $\cM(\overline{\cQ}_\Nz)$
such that $\bar{s}_a\bar{s}_b=\bar{s}_{ab}$ for $a,b \in
\Qz^\times_+$.\eprop \bproof In order to show that
$\bar{s}_n$ is unitary it suffices to show that
$1_k\bar{s}_n\bar{s}_n^*=1_k$ and
$1_k\bar{s}_n^*\bar{s}_n=1_k$ for all $k,n$.\eproof

We can also extend the generating unitary $u$ in $\cQ_\Nz$ to
a unitary in the multiplier algebra
$\cM(\overline{\cQ}_\Nz)$. We define the unitary $\bar{u}$ in
$\cM(\overline{\cQ}_\Nz)$ by the identity
$$\bar{u}\;1_k =\iota_k(u^k)$$
We can also define fractional powers of $\bar{u}$ by setting
$$\bar{u}^{1/n}\;1_{kn} =\iota_{kn}(u^k)$$

\bprop For all $a\in \Qz_+^\times$ and $b\in \Qz$ we have the
identity
$$\bar{s}_a \bar{u}^b = \bar{u}^{ab} \bar{s}_a$$
\eprop \bproof Check that
$\bar{s}_{1/n}\bar{u}=\bar{u}^{1/n}\bar{s}_{1/n}$\eproof

Following Bost-Connes we denote by $P_{\Qz}^+$ the
$ax+b$-group
$$P_{\Qz}^+ = \{ \left(\begin{array}{cc}1&b\\0&a\end{array}\right)
\mathop{|} a\in\Qz^\times_+,b\in\Qz\}$$ It follows from the
previous proposition that we have a representation of
$P_{\Qz}^+$ in the unitary group of $\cM(\bar{\cQ}_\Nz)$.
Denote by $\Az_f$ the locally compact space of finite adeles
over $\Qz$ i.e.
$$\Az_f = \{ (x_p)_{p\in \cP} \mathop{|} x_p\in\hat{\Qz}_p
\textrm{ and } x_p\in \hat{\Zz}_p \textrm{ for almost all } p\}$$
where $\cP$ is the set of primes in $\Nz$.\mn

The canonical commutative subalgebra $D$ of $\cQ_\Nz$ is by
the Gelfand transform isomorphic to $\cC(X)$ where $X$ is the
compact space $\hat{\Zz }=\prod_{p\in\cP}\hat{\Zz}_p$. It is
invariant under the endomorphisms $\varphi_n$ and we obtain an
inductive system of commutative algebras $(D,\varphi_n)$. The
inductive limit $\bar{D}$ of this system is a canonical
commutative subalgebra of $\bar{\cQ}_\Nz$. It is isomorphic
to $\cC_0(\Az_f)$. In fact the spectrum of $\bar{D}$ is the
projective limit of the system (Spec$D,\hat{\varphi}_n)$ and
$\varphi_n$ corresponds to multiplication by $n$ on $\hat{\Zz
}$.

\btheo\label{adel} The algebra $\overline{\cQ}_\Nz$ is
isomorphic to the crossed product of $\cC_0(\Az_f)$ by the
natural action of the $ax+b$-group $P_{\Qz}^+$.\etheo \bproof
Denote by $B$ the crossed product and consider the projection
$e\in \cC_0(\Az_f)\subset B$ defined by the characteristic
function of the maximal compact subgroup
$\prod_{p\in\cP}\hat{\Zz}_p\subset \Az_f$. Consider also the
multipliers $\bar{s}_n$ and $\bar{u}$ of $B$ defined as the
images of the elements $
\left(\begin{array}{cc}1&0\\0&n\end{array}\right)$ and $
\left(\begin{array}{cc}1&1\\0&1\end{array}\right)$ of
$P_\Qz^+$, respectively. Then the elements $s_n=\bar{s}_ne$
and $u=\bar{u}e$ satisfy the relations defining $\cQ_\Nz$.
Moreover, the $C^*$-algebra generated by the $s_n$ and $u$
contains the subalgebra $e\cC_0(\Az_f)$ of $B$. This shows
that this subalgebra equals $eBe$ and, by simplicity of
$\cQ_\Nz$, it follows that $eBe\cong \cQ_\Nz$.\eproof It is a
somewhat surprising fact that we get exactly the same
$C^*$-algebra, even together with its canonical action of
$\Rz$ if we replace the completion $\Az_f$ of $\Qz$ at the
finite places by the completion $\Rz$ at the infinite place.
\btheo \label{infplace} The algebra $\overline{\cQ}_\Nz$ is
isomorphic to the crossed product of $\cC_0(\Rz)$ by the
natural action of the ax+b-group $P_{\Qz}^+$. There is an
isomorphism $\overline{\cQ}_\Nz\to \cC_0(\Rz)\rtimes
P_{\Qz}^+$ which carries the natural one parameter group
$\lambda_t$ to a one parameter group $\lambda'_t$ such that
$\lambda'_t(fu_g)=a^{it}fu_g$ where $u_g$ denotes the unitary
multiplier of $\cC_0(\Rz)\rtimes P_{\Qz}^+$ associated with an
element $g\in P_{\Qz}^+$ of the form
$$g=\left(\begin{array}{cc} 1&b\\0&a \end{array}\right)$$
\etheo In order to prove this we need a little bit of
preparation concerning the representation of the
Bunce-Deddens algebra $\cF$ as an inductive limit. As noted
above it is an inductive limit of the inductive system
$(M_n(\cC (S^1)))$ with maps sending the unitary generator
$z$ of $\cC (S^1)$ to a unitary $v$ in $M_k(\Cz)\otimes\cC
(S^1)$ satisfying $v^k =1\otimes z$. We observe now that such
a $v$ is unique up to unitary equivalence.

\blemma\label{uniq} Let $v_1,v_2$ be two unitaries in
$M_k(\Cz)\otimes\cC (S^1)$ such that $v_1^k=v_2^k=1\otimes
z$. Then there is a unitary $w$ in $M_k(\Cz)\otimes\cC (S^1)$
such that $v_2=wv_1w^*$.\elemma

\bproof The spectral projections $p_i(t)$ for the different
$k-th$ roots of $e^{2\pi it}1$, given by $v_1(t)$ and
$v_2(t)$, have to be continuous functions of $t$. This implies
that all the $p_i(t)$ are one-dimensional and that, after
possibly relabelling, we must have the situation where
$$p_i(t+1) = p_{i+1}(t)$$
This means that the $p_i$ combine to define a line bundle on
the $k$-fold covering of $S^1$ by $S^1$. However any two such
line bundles are unitarily equivalent.\eproof

We now determine the crossed product $\cC_0(\Rz)\rtimes \Qz$ where
$\Qz$ acts by translation.

\blemma\label{cros} The crossed product $\cC_0(\Rz)\rtimes
\Qz$ is isomorphic to the stabilized Bunce-Deddens algebra
$\cK\otimes\cF$.\elemma \bproof The algebra $\cC_0(\Rz)\rtimes
\Qz$ is an inductive limit of algebras of the form
$\cC_0(\Rz)\rtimes \Zz$ with respect to the maps
$$\beta_k :\cC_0(\Rz)\rtimes \Zz\lori \cC_0(\Rz)\rtimes \Zz$$
obtained from the embeddings $\Zz\cong
\Zz\frac{1}{k}\hookrightarrow \Qz$. It is well known that
$\cC_0(\Rz)\rtimes \Zz$ is isomorphic to $\cK\otimes
\cC(S^1)$. An explicit isomorphism is obtained from the map
$$\sum_{n\in \Zz}f_nu^n \mapsto \sum_{n\in \Zz,k\in
\Zz}\tau_k(f_n)e_{k,k+n}$$ where $\tau_k$ denotes translation
by $k$, $u^k$ denotes the unitary in the crossed product
implementing this automorphism and $e_{ij}$ denote the matrix
units in $\cK\cong \cK (\ell^2\Zz)$.

This map sends $\cC_0(\Rz)\rtimes \Zz$ to the mapping torus algebra
$$\{f\in \cC(\Rz,\cK) \mathop{|} f(t+1)=Uf(t)U^*\}\cong
\cK\otimes \cC(S^1)$$ where $U$ is the multiplier of
$\cK=\cK(\ell^2\Zz)$ given by the bilateral shift on
$\ell^2(\Zz)$.

A projection $p$ corresponding, under this isomorphism, to
$e\otimes 1$ in $\cK\otimes \cC(S^1)$ with $e$ a projection
of rank 1 can be represented in the form $p=ug + f + gu^*$
with appropriate positive functions $f$ and $g$ with compact
support on $\Rz$. Under the map $\cC_0(\Rz)\rtimes
\Zz\cong\cC_0(\Rz)\rtimes k\Zz\lori \cC_0(\Rz)\rtimes \Zz$,
the projection $p$ is mapped to $p'=u^kg_k + f_k + g_ku^{*k}$
where $g_k(t):= g(t/k)$, $f_k(t):= f(t/k)$. Now, $p'$
corresponds to a projection of rank $k$ in $\cK\otimes
\cC(S^1)$.

Let $z$ be the unitary generator of $\cC (S^1)$. Then the
element $e\otimes z$ corresponds the function $e^{2\pi it}(ug
+ f + gu^*)$ which is mapped to $v=e^{2\pi it/k}(u^kg_k + f_k
+ g_ku^{*k})$. Thus $v^k=p'$ and $p'$ corresponds to a
projection of rank $k$ in $\cK\otimes \cC (S^1)$. \mn On the
other hand $\cF= \mathop{\lim}\limits_{{
{\textstyle\longrightarrow} }}A_n$ where $A_n=M_n(\cC(S^1))$
and the inductive limit is taken relative to the maps
$\alpha_k:M_n(\cC(S^1))\to M_{kn}(\cC(S^1))$ which map the
unitary generator $z$ of $\cC(S^1)$ to an element $v$ such
that $v^k=1\otimes z$. Compare this now to the inductive
system $A'_n= \cC_0(\Rz)\rtimes \Zz\frac {1}{n}$ with respect
to the maps $\beta_k$ considered above. From our analysis of
$\beta_k$ and from Lemma \ref{uniq} we conclude that there are
unitaries $W_k$ in $\cM( A_{kn}')$ such that the following
diagram commutes
$$\xymatrix{ A_n\quad\ar[d]\ar[r]^{\alpha_k}&
\quad A_{kn}\ar[d]\\
A_n'\quad\ar[r]^{\Ad W_k\beta_k} & \quad A'_{kn}}$$ where the
vertical arrows denote the natural inclusions
$M_n(\cC(S^1))\to \cK\otimes\cC(S^1)$.\mn We conclude that
these natural inclusions induce an isomorphism from
$\cK\otimes \cF = \mathop{\lim}\limits_{{
{\textstyle\longrightarrow} }}\cK\otimes A_n$ to
$\cC_0(\Rz)\rtimes \Qz = \mathop{\lim}\limits_{{
{\textstyle\longrightarrow} }}A_n'$.\eproof

\bproof of Theorem \ref{infplace}. Consider the commutative
diagram and the injection $\cF\hookrightarrow
\cC_0(\Rz)\rtimes \Qz$ constructed at the end of the proof of
Lemma \ref{cros}. Note also that the $\beta_k$ are
$\sigma$-unital and therefore extend to the multiplier
algebra. This shows that the injection transforms $\alpha_p$
into $\beta_p$ times an approximately inner automorphism, for
each prime $p$. Therefore the injection transforms
$\lambda_t$, which is the dual action on the crossed product
for the character $n\mapsto n^{it}$ of $\Nz^\times$, to the
restriction of the corresponding dual action on the crossed
product $(\cC_0(\Rz)\rtimes \Qz)\rtimes \Nz^\times$.

Finally, note that the endomorphisms $\beta_n$ of
$\cC_0(\Rz)\rtimes \Qz$ are in fact automorphisms and that
therefore $\beta$ extends from a semigroup action to an
action of the group $\Qz^\times_+$. This shows that
$(\cC_0(\Rz)\rtimes \Qz)\rtimes \Nz^\times =
(\cC_0(\Rz)\rtimes \Qz)\rtimes \Qz_+^\times
=\cC_0(\Rz)\rtimes P_\Qz^+$. \eproof

\bremark If we consider the full space of adeles $\Qz_\Az =
\Az_f\times\Rz $ rather than the finite adeles $\Az_f$ or the
completion at the infinite place $\Rz$, we obtain the
following situation. By \cite{Weil}, IV, \S 2, Lemma 2, $\Qz$
is discrete in $\Qz_\Az$ and the quotient is the compact space
$(\hat{\Zz}\times \Rz)/\Zz$. Thus, the crossed product
$\cC_0(\Qz_\Az)\rtimes\Qz$ which would be analoguous to the
Bunce-Deddens algebra is Morita equivalent to
$\cC((\hat{\Zz}\times \Rz)/\Zz)$. Now, $(\hat{\Zz}\times
\Rz)/\Zz$ has a measure which is invariant under the action
of the multiplicative semigroup $\Nz^\times$. Therefore the
crossed product $\cC_0(\Qz_\Az)\rtimes P_\Qz^+$ has a trace
and is not isomorphic to $\cQ_\Nz$ .\eremark

\section{The case of the multiplicative semigroup
$\Zz^\times$}

We consider now the analogous $C^*$-algebras where we replace
the multiplicative semigroup $\Nz^\times$ by the semigroup
$\Zz^\times$. This case is also important in view of possible
generalizations from $\Qz$ to more general number fields.
Thus, on $\ell^2 (\Zz)$ we consider the isometries $s_n, n\in
\Zz^\times$ and the unitaries $u^m, m\in \Zz$ defined by
$s_n(\xi_k) = \xi_{nk}$ and $u^m (\xi_k) =\xi_{k+m}$. As
above, these operators satisfy the relations
\begin{equation}\label{rel}s_n s_m =s_{nm}\qquad s_nu^m
=u^{nm}s_n \qquad \sum_{i=0}^{n-1} u^ie_nu^{-i} =1
\end{equation} The operator $s_{-1}$ plays a somewhat special
role and we therefore denote it by $f$. Then the $s_n, n\in
\Zz$ and $u$ generate the same $C^*$-algebra as the $s_n,
n\in \Nz$ together with $u$ and $f$. The element $f$ is a
selfadjoint unitary so that $f^2=1$ and we have the relations
$$fs_n =s_n f\qquad fuf = u^{-1}$$
We consider now again the universal $C^*$-algebra generated
by isometries $s_n, n\in \Zz$ and a unitary $u$ subject to
the relations (\ref{rel}). We denote this $C^*$-algebra by
$\cQ_\Zz$. We see from the discussion above that we get a
crossed product $\cQ_\Zz \cong \cQ_\Nz\rtimes \Zz/2$ where
$\Zz/2$ acts by the automorphism $\alpha$ of $\cQ_\Nz$ that
fixes the $s_n, n\in \Nz$ and $\alpha (u) =u^{-1}$.

\btheo The algebra $\cQ_\Zz$ is simple and purely infinite.
\etheo \bproof Composing the conditional expectation $G:
\cQ_\Nz \to \cD$ used in the proof of Theorem \ref{pi} with
the natural expectation $\cQ_\Zz = \cQ_\Nz\rtimes \Zz/2 \lori
\cQ_\Nz$ we obtain again a faithful expectation
$G':\cQ_\Zz\to \cD$. The rest of the proof follows exactly
the proof of \ref{pi}, using in addition the fact that the
$f_i$ in that proof can be chosen such that $f_iff_i=0$.
\eproof

Denote by $P_\Qz$ the full $ax+b$-group over $\Qz$, i.e.
$$P_{\Qz} = \{ \left(\begin{array}{cc}1&b\\0&a\end{array}\right)
\mathop{|} a\in\Qz^\times,b\in\Qz\}$$

\btheo $\cQ_\Zz$ is isomorphic to the crossed product of
$\cC_0 (\Az _f)$ or of $\cC_0 (\Rz)$ by the natural action of
$P_\Qz$.\etheo \bproof This follows from Theorems \ref{adel}
and \ref{infplace} since $P_\Qz = P_\Qz^+\rtimes
\Zz/2$.\eproof

On $\cQ_\Zz$ we can define the one-parameter group
$(\lambda_t)$ by $\lambda_t(s_n) = n^{it}s_n, n\in
\Nz^\times$. The fixed point algebra is the crossed product
$\cF\rtimes \Zz/2$ of the Bunce-Deddens algebra by $\Zz/2$.

In order to compute the $K$-groups of $\cQ_\Zz$ we first
determine the $K$-theory for $\cF ' =\cF\rtimes \Zz/2$. This
algebra is the inductive limit of the subalgebras $A_n'=
C^*(u,f,e_n)$.

\blemma\label{K0} (a) The $C^*$-algebra $C^*(u,f)$ is isomorphic to
$C^*(D)$, where $D$ is the dihedral group $D=\Zz \rtimes \Zz/2$
($\Zz/2$ acts on $\Zz$ by $a\mapsto -a$). For each $n=1,2,\ldots$,
the algebra $A'_n$ is isomorphic
to $M_n(C^*(D))$.\\
(b) The generators of $K_0(C^*(D))$ are given by the classes
of the spectral projections $f^+$ and $(uf)^+$ of $f$ and
$uf$, for the eigenvalue 1, and by $1$. We have $K_0(A_n')=\Zz
^3$ and $K_1(A_n')= 0$
for all $n$.\\
(c) Let $p$ be prime. If $p=2$, then the map $K_0(A'_n)\to
K_0 (A_{pn}')$ is described by the matrix
$$\left(\begin{array}{ccc}  2&1&0\\0&0&1\\0&0&1
\end{array}\right)$$
If $p$ is odd, then the map $K_0(A'_n)\to K_0 (A_{pn}')$ is
described by the matrix
$$\left(\begin{array}{ccc}  p&\frac{p-1}{2}&\frac{p-1}{2}
\\0&1&0\\0&0&1
\end{array}\right)$$
\elemma \bproof (a) It is clear that the universal algebra generated
by two unitaries $u,f$ satisfying $f^2=1$ and $fuf=u^*$ is
isomorphic to $C^*(D)$.

In the decomposition of $C^*(u,f,e_2)$ with respect to the
orthogonal projections $e_2$ and $ue_2u^{-1}$, the elements
$u$ and $f$ correspond to matrices
$$\left(\begin{array}{cc}0&w\\ 1&0 \end{array}\right)\quad
\rm{and}\quad \left(\begin{array}{cc}f_0&0\\ 0&f_1
\end{array}\right)$$ where $w$ is unitary and $f_0,f_1$ are
symmetries (selfadjoint unitaries).

The relations between $u$ and $f$ imply that $wf_1=f_0$ and that
$wf_1$ is a selfadjoint unitary, whence $f_1 wf_1 = w^*$. Thus
$A'_2$ is isomorphic to $M_2(C^*(w,f_1))$ and $w,f_1$ satisfy the
same relations as $u,f$.

If $p$ is an odd prime, then in the decomposition of $C^*(u,f,e_p)$
with respect to the pairwise orthogonal projections
$e_p,ue_pu^{-1},\ldots, u^{p-1}e_pu^{-(p-1)}$, $u$ and $f$
correspond to the following $p\times p$-matrices
$$\left(\begin{array}{ccccc}
0&0&\ldots&&w\\
1&0&\ldots&&0\\
0&1&\ldots&&0\\
&&\ldots&&\\
0&\ldots&&1&0 \end{array}\right)\quad \rm{and}\quad
\left(\begin{array}{ccccc}
f_0&0&\ldots&&0\\
0&0&\ldots&&f_1^*\\
0&0&\ldots&f_2^*&0\\
&&\ldots&&\\
0&f_2&\ldots &&\\
f_1&\ldots&&0&0 \end{array}\right)$$ where $w$ and the
$f_1,\ldots,f_{\frac{p-1}{2}}$ are unitary and $f_0$ is a symmetry.

The relation $fuf=u^*$ implies that $f_0 =wf_1$, $f_1=f_2=\ldots
=f_{p-1}$ and $f_i^2=1$ for all $i$. From $(wf_1)^2=1$ we derive
that $f_1wf_1=w^*$. Thus $A'_p\cong M_p(C^*(w,f_1))$ and $w, f_1$
satisfy the same relations as $u,f$.

The case of general $n$ is obtained by iteration using the fact that
$C^*(u,f,e_{nm})\cong M_m(C^*(u,f,e_{n}))$.

(b) This follows for instance from the well known fact that $C^*(D)$
$\cong (\Cz\star\Cz)^\sim$.

(c) Let $p=2$. Using the description of $K_0(C^*(u,f))$ under (b)
and the description of the map $C^*(u,f)\to C^*(u,fe_2)$ from the
proof of (b), we see that the generators $[1]$, $[(uf)^+]$ and
$[f^+]$ are respectively mapped to $[1]$, $[1]$ and
$[f^+]+[(uf)^+]$.

Let $p$ be an odd prime. The matrix corresponding to $f$ in the
proof of (b) is conjugate to the matrix where all the $f_i$,
$i=1,\ldots, (p-1)/2$  are replaced by 1. Thus the class of $f^+$ is
mapped to $[(wf_1)^+]+[1]$.

The matrix corresponding to $uf$ is conjugate to a second diagonal
matrix with all entries 1, except for the middle entry which is
$f_1$. Thus the class of $(uf)^+$ is mapped to $[f_1^+]+
\frac{p-1}{2}[1]$.

\eproof

\bprop The $K$-groups of the algebra $\cF '=\cF\rtimes\Zz/2$
are given by $K_0(\cF')=\Qz\oplus\Zz$ and
$K_1(\cF')=0$.\eprop \bproof This follows immediately from the
fact that $\cF'$ is the inductive limit of the algebras
$A_n'$ and the description of the maps $K_0(A_n')\to
K_0(A_{pn}')$ given in Lemma \ref{K0} (c).\eproof

In analogy to the case over $\Nz$ we denote by $B_n'$ the
$C^*$-subalgebra of $\cQ_\Zz$ generated by
$u,f,s_{p_1},\ldots,s_{p_n}$. We now immediately deduce

\btheo We have
$$K_0(B_n')=\Zz^{2^{n-1}}\qquad K_1(B_n')=\Zz^{2^{n-1}}$$
and
$$K_0(\cQ_\Zz)=\Zz^\infty\qquad K_1(\cQ_\Zz)=\Zz^\infty$$\etheo

\end{document}